\documentclass{amsart}

\usepackage{amssymb}
\theoremstyle{plain}

\newtheorem{thm}{Theorem}

\newtheorem{cor}[thm]{Corollary}
\newtheorem{lem}[thm]{Lemma}

\theoremstyle{definition}
\newtheorem{df}[thm]{Definition}

\newtheorem{calc}[thm]{Calculation}

\author{Michael Robinson}
\address{Center for Applied Mathematics\\657 Rhodes Hall\\Cornell
  University, Ithaca, NY 14850}
\email{robinm@cam.cornell.edu}


\title[Classification of connecting solutions of parabolic
equations]{Classification of connecting solutions of semilinear
parabolic equations}

\begin{document}

\begin{abstract}
For a given semilinear parabolic equation with polynomial
nonlinearity, many solutions blow up in finite time.  For
a certain large class of these equations, we show that some of the
solutions which do not blow up actually tend to equilibria.  The
characterizing property of such solutions is a finite energy
constraint, which comes about from the fact that this class of
equations can be written as the $L^2$ gradient of a certain
functional.
\end{abstract}

\maketitle

\section{Introduction}

In this article, the global behavior of smooth solutions to the
semilinear parabolic equation
\begin{equation}
\label{pde}
\frac{\partial u(x,t)}{\partial t}=\Delta u(x,t) - u^N + \sum_{i=0}^{N-1}
a_i(x) u^i(t,x)=\Delta u + P(u), \text{ for } (t,x) \in \mathbb{R}^{n+1}
\end{equation}
is considered, where $N\ge 2$ and $a_i \in L^1 \cap
L^\infty(\mathbb{R}^n)$ are smooth with all derivatives of all orders
bounded.  Since the linear portion of the right side of \eqref{pde} is
a sectorial operator, we can use \eqref{pde} to define a nonlinear
semigroup. \cite{Henry} \cite{Mazya_2007} Indeed, in
\cite{RobinsonIMEX}, it is shown that short-time solutions exist to
\eqref{pde} when initial conditions lie in a certain subset of $L^1
\cap L^\infty$.  This turns \eqref{pde} into a dynamical system, the
behavior of which is largely controlled by its equilibria.  In
particular, our main result is that solutions to \eqref{pde} which
connect two equilibrium solutions of \eqref{pde} in a certain strong
sense are characterized by finite energy (Definition \ref{energy_df}).
The study of this kind of problem is not entirely new.  Blow-up
behavior for equations like \eqref{pde} was examined in a classic
paper by Fujita. \cite{Fujita} For somewhat more restricted
nonlinearities, Du and Ma were able to use squeezing methods to obtain
similar results to what we obtain here.  In particular, they also show
that certain kinds of solutions approach equilibria. \cite{DuMa2001}

\section{Finite energy constraints}

It is well-known that solutions to \eqref{pde} exist along strips of
the form $(t,x)\in I\times \mathbb{R}^n$ for sufficiently small
$t$-intervals $I$.  One might hope to extend such solutions to all of
$\mathbb{R}^{n+1}$, but for certain choices of initial conditions,
such global solutions may fail to exist.  Fujita's classic paper
\cite{Fujita} gives examples of this ``blow-up'' pathology.  We will
specifically avoid it by considering only global solutions to
\eqref{pde}.

Our analysis of \eqref{pde} will make considerable use of the fact
that it is a gradient differential equation.  That is, solutions to
\eqref{pde} are integral curves for the gradient of a certain
functional in $L^2(\mathbb{R}^n)$.

\begin{df}
\label{action_df}
Observe that the right side of \eqref{pde} is the $L^2(\mathbb{R}^n)$
  gradient of the following {\it action functional}:
\begin{equation}
\label{action_eqn}
A(f)=\int \frac{1}{2} \|\nabla f\|^2 - \frac{u^{N+1}}{N+1}+ \sum_{i=0}^{N-1}
\frac{a_i(x)}{i+1} u^{i+1}(t,x) dx.
\end{equation}
\end{df}

It is then evident that along a solution $u(t)$ to \eqref{pde},
$A(u(t))$ is a monotone function.  As an immediate consequence,
nonconstant $t$-periodic solutions to \eqref{pde} do not exist.  

\begin{df}
\label{energy_df}
The {\it energy functional} is the following quantity defined on a
dense subset of $L^2(\mathbb{R}^{n+1})$:
\begin{equation}
\label{energy_eqn}
E(u)=\frac{1}{2}\int_{-\infty}^\infty \int \left | \frac{\partial u}{\partial t}
\right |^2 + \left | \Delta u + P(u) \right |^2 dx\, dt.
\end{equation}
\end{df}

\begin{calc}
\label{action_energy_calc}
Suppose $u \in L^2(\mathbb{R}^{n+1})$ is in the domain of definition
for the energy functional, then
\begin{eqnarray*}
E(u)&=&\frac{1}{2}\int_{-T}^T \int \left | \frac{\partial u}{\partial t}
\right |^2 + \left | \Delta u + P(u) \right |^2 dx\, dt\\
&=&\frac{1}{2}\int_{-T}^T \int \left ( \frac{\partial u}{\partial t}
 - \Delta u - P(u) \right )^2 + 2\frac{\partial u}{\partial t}\left( \Delta u + P(u) \right) dx\, dt\\
&=&\frac{1}{2}\int_{-T}^T \int \left ( \frac{\partial u}{\partial t}
 - \Delta u - P(u) \right )^2 dx\, dt + \int_{-T}^T \left< \frac{\partial
  u}{\partial t}, \Delta u + P(u) \right > dt\\
&=&\frac{1}{2}\int_{-T}^T \int \left ( \frac{\partial u}{\partial t}
 - \Delta u - P(u) \right )^2 dx\, dt + \int_{-T}^T \left< \frac{\partial
  u}{\partial t}, \nabla A(u(t))\right > dt\\
&=&\frac{1}{2}\int_{-T}^T \int \left ( \frac{\partial u}{\partial t}
 - \Delta u - P(u) \right )^2 dx\, dt + \int_{-T}^T
\frac{d}{dt} A(u(t)) dt\\
&=&\frac{1}{2}\int_{-T}^T \int \left ( \frac{\partial u}{\partial t}
 - \Delta u - P(u) \right )^2 dx\, dt + A(u(T))-A(u(-T)).\\
\end{eqnarray*}
This calculation shows that finite energy solutions to \eqref{pde}
minimize the energy functional.  If a solution to \eqref{pde} connects
two equilibria, then the energy functional measures the difference between
the values of the action functional evaluated at the two equilibria.
The main result of this article is to show the converse, so that
finite energy characterizes the solutions which connect equilibria.
\end{calc}

It is well-known that when equations like \eqref{pde} exhibit the
correct symmetry, they can support travelling wave
solutions. \cite{FiedlerScheel} A typical travelling wave solution $u$
has a symmetry like $u(t,x)=U(x-ct)$ for some $c \in \mathbb{R}$.  As
a result, it is immediate that travelling waves will have infinite
energy.  On the other hand, they also evidently connect equilibria.
As a result, Calculation \ref{action_energy_calc} shows that a
necessary condition for travelling waves is that there exists at least
one equilibrium whose action is infinite.  In this article, we will
consider only equilibria with finite action, and solutions with finite
energy.  As a result, we will not be working with travelling waves.

\section{Convergence to equilibria}

In this section, we show that finite energy solutions tend to
equilibria as $|t| \to \infty$.  In doing this, we follow Floer in
\cite{Floer_gradient} which leads us through an essentially standard
parabolic bootstrapping argument.

\begin{lem}
\label{polybound_lem}
Let $U \subseteq \mathbb{R}^n$ and $u \in W^{k,p}(U)$ satisfy $\|D^j
u\|_\infty \le C < \infty$ for $0 \le j \le k$ (in particular, $u$ is
bounded).  If $P(u) = \sum_{i=1}^N a_i u^i$ with $a_i \in L^\infty(U)$
then there exists a $C'$ such that $\|P(u)\|_{k,p} \le C' \|u\|_{k,p}$.  
\begin{proof}
First, using the definition of the Sobolev norm,
\begin{equation*}
\|P(u)\|_{k,p}=\sum_{j=0}^k \|D^j P(u) \|_p \le \sum_{j=0}^k
\sum_{i=1}^N \|D^j a_i u^i \|_p.
\end{equation*}
Now $|D^j a_i u^i| \le P_{i,j}(u,Du,...,D^j u)$ is a polynomial in $j$
variables with constant coefficients, which has no constant term.
(The constant coefficients is a consequence of the bounded derivatives
of the $a_i$.)  Additionally,
\begin{eqnarray*}
\|(D^m u)^q D^j u\|_p &=& \left ( \int \left | (D^m u)^q D^j u \right
|^p \right)^{1/p}\\
&\le& \|D^m u\|_\infty^q \left ( \int \left | D^j u \right
|^p \right)^{1/p} \le C^q \|D^j u\|_p,\\
\end{eqnarray*}
so by collecting terms,
\begin{equation*}
\|P(u)\|_{k,p} \le  \sum_{j=0}^k
\sum_{i=1}^N \|D^j a_i u^i \|_p \le \sum_{j=0}^k A_j \|D^j u\|_p \le C' \|u\|_{k,p}.
\end{equation*}
\end{proof}
\end{lem}

The following result is a parabolic bootstrapping argument that does
most of the work.  In it, we follow Floer in \cite{Floer_gradient},
replacing ``elliptic'' with ``parabolic'' as necessary.

\begin{lem}
\label{parabolic_bootstrap_lem}
If $u$ is a finite energy solution to \eqref{pde} with $\|D^j
u\|_{L^\infty((-\infty,\infty)\times V)} \le C < \infty$ for $0 \le j
\le k$ with $k\ge 1$ on each compact $V \subset \mathbb{R}^n$, then
each of $\lim_{t \to \pm \infty} u(t,x)$ exists, and converges with
$k$ of its first derivatives uniformly on compact subsets of
$\mathbb{R}^n$.  Further, the limits are equilibrium solutions to
\eqref{pde}.
\begin{proof}
Define $u_m(t,x)=u(t+m,x)$ for $m=0,1,2...$.  Suppose $U \subset
\mathbb{R}^{n+1}$ is a bounded open set and $K \subset U$ is compact.  Let
$\beta$ be a bump function whose support is in $U$ and takes the value
1 on $K$.  We take $p>1$ such that $kp>n+1$.  Then we can consider
$u_m \in W^{k,p}(U)$ (recall that $u$ and its first $k$ derivatives of
$u$ are bounded on the closure of $U$), and we have
\begin{equation*}
\|u_m\|_{W^{k+1,p}(K)} \le \|\beta u_m \|_{W^{k+1,p}(U)}.
\end{equation*}
Then using the standard parabolic regularity for the heat operator,
\begin{equation*}
\|\beta u_m \|_{W^{k+1,p}(U)} \le C_1 \left \| \left (
\frac{\partial}{\partial t} - \Delta\right ) (\beta u_m)\right \|_{W^{k,p}(U)}.
\end{equation*}
Let $P'(u)=-u^N+\sum_{i=1}^{N-1} a_i u^i$, noting carefully that we have left
out the $a_0$ term.  The usual product rule, and a little
work, as suggested in \cite{Salamon_1990} yields the following sequence of
inequalities
\begin{eqnarray*}
\|u_m\|_{W^{k+1,p}(K)}&\le& C_1 \left \|\beta
\left(\frac{\partial}{\partial t} - \Delta \right) u_m \right
\|_{W^{k,p}(U)} + C_2 \|u_m\|_{W^{k,p}(U)}\\
&\le& C_1 \left \|\beta
\left(\frac{\partial}{\partial t} - \Delta \right) u_m + \beta P'(u_m)
- \beta P'(u_m) \right
\|_{W^{k,p}(U)} + C_2 \|u_m\|_{W^{k,p}(U)}\\
&\le& C_1  \|\beta a_0\|_{W^{k,p}(U)} + C_1 \| \beta P'(u_m)
\|_{W^{k,p}(U)} + C_2 \|u_m\|_{W^{k,p}(U)}\\
&\le& C_1 \|\beta a_0\|_{W^{k,p}(U)} + C_3 \|u_m\|_{W^{k,p}(U)},\\
\end{eqnarray*}
where the last inequality is a consequence of Lemma
\ref{polybound_lem}.  By the hypotheses on $u$ and $a_0$, this implies
that there is a finite bound on $\|u_m\|_{W^{k+1,p}(K)}$, which is
independent of $m$.  Now by our choice of $p$, the general
Sobolev inequalities imply that $\|u_m\|_{C^{k+1-(n+1)/p}(K)}$ is
uniformly bounded.  By choosing $p$ large enough, there is a
subsequence $\{v_{m'}\} \subset \{u_m\}$ such that $v_{m'}$ and its
first $k$ derivatives converge uniformly on $K$, say to $v$.  For any
$T>0$, we observe
obtain
\begin{eqnarray*}
\int_{-T}^T \int \left | \frac{\partial v}{\partial t} \right |^2 dx
\, dt &=& \lim_{{m'} \to \infty}\int_{-T}^T \int \left | \frac{\partial v_{m'}}{\partial t} \right |^2 dx
\, dt\\
&=&\lim_{{m'} \to \infty} \int_{{m'}-T}^{{m'}+T} \int \left | \frac{\partial u}{\partial t} \right |^2 dx
\, dt=0,\\
\end{eqnarray*}
where the last equality is by the finite energy condition.  Hence
$\left|\frac{\partial v}{\partial t}\right| = 0$ almost everywhere,
which implies that $v$ is an equilibrium and that $\lim_{t\to\infty}
u(t,x) = v(x)$.  Similar reasoning works for
$t\to - \infty$.
\end{proof}
\end{lem}

Now we would like to relax the bounds on $u$ and its derivatives, by
showing that they are in fact consequences of the finite energy
condition.

\begin{lem}
\label{finite_energy_consequences_lem}
Suppose that either $n=1$ (one spatial dimension) or $N$ is odd, then
we have the following. If $u$ is a finite energy solution to
\eqref{pde}, then for every $v=(t_0,x_0)\in \mathbb{R}^{n+1}$, the
limits $\lim_{s\to\pm\infty} u(t+st_0,x+s x_0)$ exist
uniformly on compact subsets of $\mathbb{R}^{n+1}$, and
additionally,
\begin{itemize}
\item $u$ is bounded,
\item the derivatives $Du$ are bounded,
\item and therefore the limits are continuous equilibrium solutions.
\end{itemize}
\begin{proof}
Note that since 
\begin{equation*}
E(u)=\frac{1}{2}\int_{-\infty}^\infty \int \left | \frac{\partial u}{\partial t}
\right |^2 + \left | \Delta u + P(u) \right |^2 dx\, dt < \infty,
\end{equation*}
we have that for any $\epsilon>0$,
\begin{equation*}
\lim_{T\to\infty}\frac{1}{2}\int_{T-\epsilon}^{T+\epsilon} \int \left | \frac{\partial u}{\partial t}
\right |^2 + \left | \Delta u + P(u) \right |^2 dx\, dt = 0,
\end{equation*}
whence $\lim_{t\to\infty} \left|\frac{\partial u}{\partial t}\right| =
0$ for almost all $x$.  So this gives that the limit is an equilibrium
almost everywhere.  Of course, this argument works for $t\to -\infty$.

Now in the case of $N$ being odd, a comparison principle shows that
solutions to \eqref{pde} are always bounded.  So we need to consider
the case with $N$ even.  In that case, a comparison principle on
\eqref{pde} shows that $u$ is bounded from {\it above}. On the other
hand, if $N$ is even we have assumed that $n=1$ in this case, and it
follows from an easy ODE phase-plane argument that unbounded
equilibria are bounded from {\it below}.  (Here we have used that the
coefficients $a_i$ are bounded.)  As a result, we must conclude that
if a solution to \eqref{pde} tends to any equilibrium, that
equilibrium (and hence $u$ also) must be bounded.

Now observe that $\left|\frac{\partial u}{\partial t}\right|\to 0$ as
$t\to\infty$ on almost all of any compact $K \subset \mathbb{R}^n$,
and that $\left|\frac{\partial u}{\partial t}\right| \le a < \infty$
for some finite $a$ on $\{(t,x)|t=0,x\in K\}$ by the smoothness of
$u$.  By the compactness of $K$, this means that if
$\left\|\frac{\partial u}{\partial
    t}\right\|_{L^\infty((-\infty,\infty)\times K)} = \infty$, there
must be a $(t^*,x^*)$ such that $\lim_{(t,x)\to
  (t^*,x^*)}\left|\frac{\partial u}{\partial t}\right| = \infty$.
This contradicts smoothness of $u$, so we conclude
$\left|\frac{\partial u}{\partial t}\right|$ is bounded on the strip
$(-\infty,\infty)\times K$.  On the other hand, the finite energy
condition also implies that for each $v \in \mathbb{R}^n$,
\begin{equation*}
\lim_{s\to\infty} \int_{-\infty}^\infty \int_{K+sv}
\left|\frac{\partial u}{\partial t}\right| dx\, dt = 0,
\end{equation*}
whence we must conclude that $\lim_{s\to\infty}\left|\frac{\partial
    u(t,x+sv)}{\partial t}\right|= 0$ for almost every
$t\in\mathbb{R}$ and $x\in K$.  Thus the smoothness of $u$ implies
that $\left|\frac{\partial u}{\partial t}\right|$ is bounded on all of
$\mathbb{R}^{n+1}$.  

Next, note that since $\left|\frac{\partial u}{\partial t}\right|$ and
$u$ are both bounded, then so is $\Delta u$.  (Use the boundeness of
the coefficients of $P$.)  Taken together, this implies that all the
spatial first derivatives of $u$ are also bounded.

As a result, we have on $K$ a bounded equicontinuous family of
functions, so Ascoli's theorem implies that they (after extracting a
suitable subsequence) converge uniformly on compact subsets of $K$ to
a continuous limit.
\end{proof}
\end{lem}

\begin{cor}
\label{limits_to_equilibria}
Suppose that either $n=1$ or $N$ is odd.  A smooth global solution $u$
to \eqref{pde} has finite energy if and only if for any
$v=(t_0,x_0)\in \mathbb{R}^{n+1}$, each of $\lim_{s \to \pm \infty}
u(t+st_0,x+sx_0)$ exists, and converges with its first derivatives
uniformly on compact subsets of $\mathbb{R}^{n+1}$ to bounded,
continuous, finite action equilibrium solutions to \eqref{pde}.
\end{cor}

\section{Discussion}

The point of employing the bootstrapping argument of Lemma
\ref{parabolic_bootstrap_lem} is only to extract uniform convergence
of the derivatives of the solution.  As can be seen from the proof of
Lemma \ref{finite_energy_consequences_lem}, such regularity arguments
are unneeded to obtain good convergence of the solution only.

While Corollary \ref{limits_to_equilibria} is probably true for
all spatial dimensions, the proof given here cannot be generalized to
higher dimensions.  In particular, V\'{e}ron in \cite{Veron_1996}
shows that in the case of $P(u)=-u^N$, there are solutions to the
equilibrium equation $\Delta u - u^N=0$ which are {\it unbounded
below} and {\it bounded above} when the spatial dimension is greater
than one.  This breaks the proof of Lemma
\ref{finite_energy_consequences_lem}, that the limiting equilibria of
finite energy solutions are bounded for $N$ even, since the proof
requires exactly the opposite.

On the other hand, the case of $P(u)=-u|u|^{N-1}+ \sum_{i=0}^{N-1}a_i
u^i$ is considerably easier than what we have considered here.  In
particular, all solutions to \eqref{pde} are then bounded.  In that
case, the proof of Lemma \ref{finite_energy_consequences_lem} works
for all spatial dimensions.

\bibliography{classify_bib}

\begin{thebibliography}{1}

\bibitem{DuMa2001}
Yihong Du and Li~Ma.
\newblock Logistic type equations on $\mathbb{R}^n$ by a squeezing method
  involving boundary blow-up solutions.
\newblock {\em J. London Math. Soc.}, 2(64):107--124, 2001.

\bibitem{FiedlerScheel}
Bernold Fiedler and Arnd Scheel.
\newblock Spatio-temporal dynamics of reaction-diffusion equations.
\newblock In M.~Kirkilionis, R.~Rannacher, and F.~Tomi, editors, {\em Trends in
  Nonlinear Analysis}, pages 23--152. Springer-Verlag, Heidelberg, 2003.

\bibitem{Floer_gradient}
Andreas Floer.
\newblock The unregularized gradient flow of the symplectic action.
\newblock {\em Comm. Pure Appl. Math.}, 41:775--813, 1988.

\bibitem{Fujita}
Hiroshi Fujita.
\newblock {On the blowing up of solutions of the Cauchy problem for $u_t=\Delta
  u + u^{1+\alpha}$}.
\newblock {\em Tokyo University Faculty of Science Journal}, 13:109--124,
  December 1966.

\bibitem{Henry}
Dan Henry.
\newblock {\em Geometric Theory of Semilinear Parabolic Equations}.
\newblock Springer-Verlag, New York, 1981.

\bibitem{Mazya_2007}
Vladimir Maz'ya.
\newblock {Analytic criteria in the qualitative spectral analysis of the
  Schroedinger operator, {\tt arXiv:math.SP/0702427}}.

\bibitem{RobinsonIMEX}
Michael Robinson.
\newblock {IMEX} method convergence for a semilinear parabolic equation.
\newblock {\em J. Differential Equations}, 241(2):225--236, 2007.

\bibitem{Salamon_1990}
Dietmar Salamon.
\newblock Morse theory, the {C}onley index, and {F}loer homology.
\newblock {\em Bull. London Math. Soc.}, 22:113--140, 1990.

\bibitem{Veron_1996}
Laurent V\'{e}ron.
\newblock {\em Singularities of solutions of second order quasilinear
  equations}.
\newblock Addison Wesley Longman, Essex, 1996.

\end{thebibliography}
\bibliographystyle{plain}

\end{document}